\documentclass[12pt]{article}
\usepackage{amssymb}
\usepackage{amsmath}
\usepackage{amsbsy}
\newtheorem{lemma}{Lemma}[section]

\newtheorem{theorem}{Theorem}[section]
\newtheorem{corollary}{Corollary}[section]

\def\proclaim#1{\par \bigskip\noindent {\bf #1}\bgroup\it\ }
\def\endproclaim{\egroup\par\bigskip}

\newbox\TempBox \newbox\TempBoxA

\def\pr{\textsf{P}} 
\def\ep{\textsf{E}} 
\def\Cov{\textsf{Cov}} 
\def\Var{\textsf{Var}} 

\def\text#1{\mbox{\rm #1}}

\def\underwiggle 1{
\ifmmode\setbox\TempBox=\hbox{$ 1$}\else\setbox\TempBox=\hbox{
1}\fi \setbox\TempBoxA=\hbox to \wd\TempBox{\hss\char'176\hss}
\rlap{\copy\TempBox}\smash{\lower9pt\hbox{\copy\TempBoxA}} }

\begin{document}

\begin{center}
{\Large Cram\'{e}r Type Moderate Deviation for the Maximum of the Periodogram with Application to Simultaneous Tests in Gene Expression Time
Series}
\end{center}

\begin{center}
{{Weidong Liu\footnote{Email: liuweidong99@gmail.com} and Qi-Man Shao\footnote{Email: maqmshao@ust.hk. Research partially supported by Hong Kong
RGC CERG 220608}  }}

{{\it Hong  Kong University of Science and Technology } }

\end{center}

\noindent {\rm {\bf Abstract.} In this paper, Cram\'{e}r type moderate deviations for the maximum of the periodogram and its studentized version
are derived. The results are then applied  to a simultaneous testing problem in gene expression time series. It is shown that the level of the
simultaneous tests is accurate provided that the number of genes $G$ and the sample size $n$ satisfy $G=\exp(o(n^{1/3}))$.


 }
\section{Introduction}
 \setcounter{equation}{0}

Let $X_{1}, X_{2}, \cdots$ be a sequence of random variables. Define
 the periodogram ordinates for $\{X_{n}\}$ at the standard frequencies $
 \omega_{j} =2\pi j/n$ by
\begin{eqnarray*}
I_{n}(\omega_{j})=\frac{1}{n}\Big{|}\sum_{k=1}^{n}X_{k}e^{ik\omega_{j}}\Big{|}^{2},
\end{eqnarray*}
where
 $1\leq j\leq q$ and $q=[(n-1)/2].
$

The periodogram is a fundamental tool in spectral analysis and is often used to detect periodic patterns in various real applications, such as the
analysis of gene expression data and the study of earthquake. Theoretical properties of the periodogram have been extensively studied.  An, Chen
and Hannan (1983) obtained the logarithm law for the maximum of the periodogram; Davis and Mikosch (1999), Mikosch, Resnick and Samorodnitsky
(2000), Lin and Liu (2009a) obtained the asymptotic distribution for the maximum of the periodogram under the i.i.d. and linear process cases, the
heavy-tailed case and nonlinear time series case respectively; Fay and Soulier (2001) obtained central limit theorems for functionals of the
periodogram; Shao and Wu (2007) obtained asymptotic distributions for the periodogram and empirical distribution function of the periodogram for a
wide class of nonlinear processes. When $\{X_{n}\}$ are independent and identically distributed (i.i.d.) random variables with
$\Var(X_{1})=\sigma^{2}$ and $\ep |X_{1}|^{2+\delta}<\infty$ for some $\delta>0$,  Davis and Mikosch (1999) shows that
\begin{eqnarray}\label{eq:f}
\lim_{n\rightarrow\infty}\pr\Big{(}\max_{1\leq j\leq q}I_{n}(\omega_{j})/\sigma^{2}-\log q\leq y\Big{)}=\exp(-\exp(-y)).
\end{eqnarray}

The main purpose of this paper is to study the Cram\'{e}r type moderate deviations  for the maximum of the periodogram and its studentized
version. That is, what is the largest possible $a_{n}$ so that
\begin{eqnarray}\label{eq:g}
\frac{\pr\Big{(}\max_{1\leq j\leq q}I_{n}(\omega_{j})/\sigma^{2}-\log q\geq y\Big{)}}{1-\exp(-\exp(-y))}\rightarrow 1
\end{eqnarray}
uniformly in  $y\in[-\log q, a_{n}]$, or for the studentized periodogram, what is the largest possible $b_{n}$ so that
\begin{eqnarray}\label{eq:g1}
\frac{\pr\Big{(}\frac{\max_{1\leq j\leq q}I_{n}(\omega_{j})}{q^{-1}\sum_{j=1}^{q}I_{n}(\omega_{j})}-\log q\geq
y\Big{)}}{1-\exp(-\exp(-y))}\rightarrow 1
\end{eqnarray}
uniformly in  $y\in[-\log q, b_{n}]$. We shall show that $a_{n}$ depends on the moment condition of $\{X_{n}\}$. For example, if
$\ep|X_{1}|^{2+\delta}<\infty$, $\delta>0$, the largest possible value of $a_{n}$ is $\frac{\delta}{2}\log n$, but $a_{n}$ can be chosen
$o(n^{1/3})$ if the moment generating function of $X_{1}$ is finite. However, the situation becomes totally different for the studentized
periodogram. We shall prove that $b_{n}=o(n^{1/3})$ provided that $\ep X^{4}_{1}<\infty$.

 The paper is organized as follows. Our main results, Theorems 2.1-2.3, are stated in Section
2, while proofs of the main results are postponed to Section 4. Our moderate deviation results are motivated by simultaneous tests in gene
expression time series. Theoretical results for the simultaneous tests and simulation study are discussed in Section 3.

\section{Main results}
\setcounter{equation}{0}

Throughout this paper, we assume $\{X_{n}\}$ are i.i.d. random variables. Our first result is the moderate deviation for the maximum of the
periodogram for
 $y\leq c\log n$ for some $c>0$.
 Such type of moderate deviation for the partial sums of $\{X_{n}\}$ has been studied in literature, e.g. by   Michel (1976),
  Amosova (1982), Petrov (2002) and Wu and Zhao (2008).
\begin{theorem} (i)
Suppose that  for some $c>0$,
\begin{eqnarray}\label{th1}
n^{c+1}\pr\Big{(}|X_{1}|\geq \sqrt{n\log n}\Big{)}=o(1)
\end{eqnarray}
as $n\rightarrow\infty$. Then we have
\begin{eqnarray}\label{th2}
\lim_{n\rightarrow\infty}\frac{\pr\Big{(}\max_{1\leq j\leq q}I_{n}(\omega_{j})/\sigma^{2}-\log q\geq y\Big{)}}{ 1-\exp(-\exp(-y))}=1
\end{eqnarray}
uniformly in $y\in[-\log q, ~c\log n]$, where $\sigma^{2}=\Var(X_{1})$.

(ii) If for some $\sigma>0$ (\ref{th2}) holds  uniformly in $y\in[-\log q, ~c\log n]$ with some $c>0$, then we have
\begin{eqnarray}\label{th3}
n^{c+1}\pr\Big{(}|X_{1}|\geq \sqrt{n\log n}\Big{)}=O(1).
\end{eqnarray}
\end{theorem}

Theorem 2.1 (ii) shows that  condition (\ref{th1}) is nearly optimal and hence the range depends on the moment assumption. On the other hand, when
the moment generating function exists, the range can be extended to $o(n^{1/3})$.
\begin{theorem}
Assume $\Var(X_{1})=\sigma^{2}$ and $\ep e^{t_{0}|X_{1}|}<\infty$ for some $t_{0}>0$. Then
\begin{eqnarray*}
\lim_{n\rightarrow\infty}\frac{\pr\Big{(}\max_{1\leq j\leq q}I_{n}(\omega_{j})/\sigma^{2}-\log q\geq y\Big{)}}{1-\exp(-\exp(-y))}=1
\end{eqnarray*}
uniformly in  $y\in [-\log q, ~o(n^{1/3})).$
\end{theorem}

We next consider  the maximum of the studentized periodogram. Theorem 2.3 below shows that the moment conditions in Theorems 2.1 and 2.2 can be
significally reduced for the studentized version.
\begin{theorem} If $\ep X^{4}_{1}<\infty$, then
\begin{eqnarray}\label{eq:g2}
\lim_{n\rightarrow\infty}\frac{\pr\Big{(}\frac{\max_{1\leq j\leq q}I_{n}(\omega_{j})}{q^{-1}\sum_{j=1}^{q}I_{n}(\omega_{j})}-\log q\geq
y\Big{)}}{1-\exp(-\exp(-y))}=1
\end{eqnarray}
uniformly in  $y\in [-\log q, ~o(n^{1/3})).$
\end{theorem}

Since the variance $\sigma^{2}$ of $X_{1}$ is typically unknown, what used in practice is actually the studentized periodogram. So the result in
Theorem 2.3 is more appealing and useful than Theorems 2.1 and 2.2. Theorem 2.3 also shares similar properties with self-normalized partial sums
of independent random variables, which usually requires much less moment assumptions; see, Shao (1997,1999) for self-normalized large deviation
without any moment assumption and Cram\'{e}r moderate deviation under finite third moment and de la Pe\"{n}a, Lai and Shao (2009) for recent
developments in the area of self-normalized limit theory. In view of the moderate deviation for self-normalized partial sums (Shao (1999)), we
conjecture that Theorem 2.3 remains true if $\ep |X_{1}|^{3}<\infty$.

\section{Application to simultaneous tests}
\setcounter{equation}{0}

\subsection{Theoretical results}
Periodic phenomena are widely studied in biology. Recently, there are quite a lot of interests in detecting periodic patterns in gene expression
time series; see Wichert, Fokianos and Strimer (2004), Ahdesm\"{a}ki et al. (2005), Chen (2005), Glynn, Chen and Mushegian (2006) and the
references therein. Due to modern technology such as microarray experiments, the data  are usually high-dimensional and we often need to make many
statistical inference simultaneously. Let $Y_{t,g}$ denote the observed expression level of gene $g$ at time $t$, $1\leq g\leq G$ and $1\leq t\leq
n$, where $G$ is the number of genes. The sample size $n$ is usually much smaller than the number of genes. Consider the following model of
periodic gene expression
\begin{eqnarray}\label{eqbb}
Y_{t,g}=\mu_{g}+\beta_{g}\cos(\omega t +\phi)+\varepsilon_{t,g},
\end{eqnarray}
where $\beta_{g}\geq 0$, $\omega\in(0,\pi)$, $\phi\in(-\pi,\pi]$, $\mu_{g}$ is the mean expression level. For each $g$,
$\varepsilon_{1,g},\cdots,\varepsilon_{n,g}$ are i.i.d. noise sequence with  mean zero. We wish to test the null hypothesis  $H_{0,g}:$
$\beta_{g}=0$ against the alternative hypothesis $H_{1,g}:$ $\beta_{g}\neq 0$. If $H_{0,g}$ is rejected, then we identify gene $g$ with a periodic
pattern in its expression. Periodogram is often used to detect periodically expressed gene. Let
\begin{eqnarray}\label{eq0}
I^{(g)}_{n}(\omega_{j})=\frac{1}{n}\Big{|}\sum_{k=1}^{n}Y_{k,g}e^{ik\omega_{j}}\Big{|}^{2},
\end{eqnarray}
where $\omega_{j}=2\pi j/n$, $1\leq j\leq [(n-1)/2]$. Define the $g$-statistic
\begin{eqnarray*}
f_{g}=\frac{\max_{1\leq j\leq q}I^{(g)}_{n}(\omega_{j})}{\sum_{j=1}^{q}I^{(g)}_{n}(\omega_{j})},
\end{eqnarray*}
 and its null distribution  $F_{n,g}(x)=\pr(f_{g}\leq x|H_{0,g})$, where $q=[(n-1)/2]$. Under the null hypothesis and the assumption that
$\varepsilon_{1,g},\cdots,\varepsilon_{n,g}$ are i.i.d normal random variables, the exact distribution for $f_{g}$ can be found in Fisher (1929):
\begin{eqnarray}\label{eq1}
\pr(f_{g}>x|H_{0,g})=\sum_{j=1}^{p}(-1)^{j-1}C^{j}_{q}(1-jx)^{q-1}=:f_{n}(x),
\end{eqnarray}
where $p=[1/x]$. Using (\ref{eq1}),  Wichert, Fokianos and Strimer (2004) proposed the following method to identify periodically expressed genes:
 1. For each time series calculate Fisher's statistic
$f_{g}$.
 2. For each of the test statistic calculate the
corresponding $p$ value $P_{g}=f_{n}(f_{g})$.
 3. Use the method of Benjamini and Hochberg (1995) to
control the False Discovery Rate (FDR) at $\theta$. Consider the set of ordered $p$ values $P_{(1)},\cdots, P_{(G)}$ and let
\begin{eqnarray}\label{eq2}
 i_{\theta}=\max\{i: P_{(i)}\leq i\theta/G\}.
\end{eqnarray}
 Reject the null hypothesis for the time
series indexed by $\mathcal{S}=\{i: P_{i}\leq P_{(i_{\theta})}\}$.

In many applications such as those arising from bioinformatics, the noise can be remarkably non-Gaussian (Ahdesm\"{a}ki et al. (2005)). Then the
values $P_{1},\cdots, P_{G}$ are only the estimators of the true $p$-values. It is natural to ask:
\begin{quote}
{\normalsize How large $G$ can be before the accuracy of simultaneous statistical inference becomes poor?}
\end{quote}
 Similar problems have been studied in Fan, Hall and Yao (2007), where they consider
$Y_{t,g}=\mu_{g}+\varepsilon_{t,g}$, the model in (\ref{eqbb}) without periodic parts, and focused on testing $H^{'}_{0,g}: \mu_{g}=0$. Let the
true $p$ value be $P^{true}_{g}=(1-F_{n,g}(f_{g}))$. From (\ref{eq2}), $P^{true}_{(i_{\theta})}$ may be of the order $O(1/G)$. Hence $\max_{1\leq
g\leq G}|P^{true}_{g}-P_{g}|=o(1)$, implied by (\ref{eq:f}), is not enough. The required accuracy between the estimated $p$ value and the true $p$
value is
\begin{eqnarray}\label{eq4}
|P_{g}-P^{true}_{g}|I\{\mathcal{H}_{g}\}=o(P^{true}_{g}) \quad\mbox{uniformly in $1\leq g\leq G$,}
\end{eqnarray}
i.e.
\begin{eqnarray*}
\max_{1\leq g\leq G}\Big{|}\frac{P_{g}}{P^{true}_{g}}-1\Big{|}I\{\mathcal{H}_{g}\}=o(1),
\end{eqnarray*}
where $\mathcal{H}_{g}=\{P_{g}>\theta/(2G) ~or~ P^{true}_{g}>\theta/(2G)\}$. (On $\mathcal{H}^{c}_{g}$, the gene $g$ is always rejected.) Some
similar requirements as (\ref{eq4}) on simultaneous tests have been proposed by Fan, Hall and Yao (2007) and Kosorok and Ma (2007), p.1460.

Recall that  $F_{n,g}(x)=\pr(f_{g}\leq x|H_{0,g})$. By examining the proof of Theorem 2.3, we can get the following corollary.

\begin{corollary}\label{c1}
Suppose that $\min_{1\leq g\leq G}\Var(\varepsilon_{1,g})\geq \kappa$ for some $\kappa>0$ which does not depend on $G$. Further assume that
$\max_{1\leq g\leq G}\ep\varepsilon_{1,g}^{4}=O(1)$. Then the null distribution $F_{n,g}(x)$ satisfies
\begin{eqnarray*}
\max_{1\leq g\leq G}\Big{|}\frac{1-F_{n,g}((y+\log q)/q)}{1-\exp(-\exp(-y))}-1\Big{|}=o(1)
\end{eqnarray*}
uniformly in  $y\in [-\log q, ~o(n^{1/3})).$
\end{corollary}

The following lemma shows that we can replace $1-\exp(-\exp(-y))$ by $f_{n}((y+\log q)/q)$.
\begin{lemma} We have
\begin{eqnarray*}
\lim_{n\rightarrow\infty}\Big{|}\frac{f_{n}((y+\log q)/q)}{1-\exp(-\exp(-y))}-1\Big{|}=0,
\end{eqnarray*}
uniformly in  $y\in [-\log q, ~o(n^{1/3})).$
\end{lemma}

  It follows from  Corollary 3.1 and Lemma 3.1 that

\begin{theorem} Suppose the conditions in Corollary \ref{c1} are satisfied and  $G=\exp(o(n^{1/3}))$. Then
 (\ref{eq4}) holds.
\end{theorem}

Theorem 3.1 shows that the level of the simultaneous tests is accurate provided that  $G=\exp(o(n^{1/3}))$, which seems to be the correct order of
asymptotics for microarray experiments with a moderate number of samples.

Using the bootstrap and a refined expansion of $t$-statistic in Theorem 1.2 of Wang (2005), Fan, Hall and Yao (2007) shows that $\exp(o(n^{1/3}))$
can be replaced by $\exp(o(n^{1/2}))$ for the tests $H_{0,g}: \mu_{g}=0$. It would be interesting to investigate whether a similar expansion as
Theorem 1.2 of Wang (2005) holds for the   maximum of the studentized periodogram.

\subsection{Simulation study} In this section, we carry out a simple simulation study to assess the finite sample performance. We generate
$2000$ genes with 100 periodic genes for different sample sizes $n$. Consider
\begin{eqnarray*}
&&Y_{t,g}=\beta(\cos(\omega^{(g)} t)+\sin(\omega^{(g)} t))+\varepsilon_{t,g},
\quad 1\leq t\leq n, \quad 1\leq g\leq 100, \\
&& Y_{t,g}=\varepsilon_{t,g},\quad\quad\quad\quad \quad\quad\quad\qquad\qquad1\leq t\leq n,\quad101\leq g\leq 2000.
\end{eqnarray*}
 $\varepsilon_{t,g}$ will be taken as $N(0,1)$, $(\sqrt{3/5})\times t(5)$,  EXP(1), $2^{-1}\times\chi^{2}(2)$,
  where $t(5)$ has the $t$
distribution with freedom 5, EXP(1) is the exponential random variable with parameter $\mu=1$, $\chi^{2}(2)$ is chi square random variable with
freedom 2. (The constants on the left hand side of random variables are chosen so that $\Var(\varepsilon_{t,g})=1$). The FDR level $\theta$ is
chosen as $0.15$ and  $0.05$. The simulation results are based on 100 replicates.

 We only give the simulation study when
$\omega^{(g)}$ is of the form of $\omega_{i}$ for some $1\leq i\leq q$. To do this, we let $\omega^{(g)}=2\pi/10$, $n=20,50$ and $\beta=1$. The
results are summarized in Table 1, where Tot.=total count identified using FDR; Pos.= the number of true positives identified using FDR; Z=the
number of true periodic genes among the smallest 100 $p$-values  genes.  We note that when the tails of $\varepsilon_{t,g}$ are heavier than that
of Gaussian random variable, the empirical FDR (EFDR) are  lower than the target FDR,  while most of periodic genes can still be found. There are
no significant differences between Gaussian noise and other noises when $n$ is  large moderately ($n=50$). Powers increase as $n$ increases.
Overall Fisher's statistic
 is relatively robust to the noise, as indicated by our theorem.
We refer to Wichert, Fokianos and Strimer (2004) for some real data analysis.

\begin{table}[!h]
\scriptsize  {\caption{}
\begin{center}
\begin{tabular}{|r|c|r|r|r|r|r|}
\hline
  $\theta$&   &        Tot.(Pos.)  & $EFDR$ & Tot.(Pos.) & $EFDR$ \\
  \hline
 && $n=20$ && $n=50$&\\
  \hline
   0.15    &  Normal                   &39.9(33.7)~~   &0.155 ~~&116.2(99.6) &0.143 ~\\
   0.05        &                 &12.6(11.9)~~~&     0.059~~ &102.6(99.2)&0.033~ \\
  \hline
   Z        &                 &62~~~& &97.1~~~~~&\\
   \hline\hline
   0.15    &      EXP(1)                  &55.3(49.5)~~   &0.105 ~~&104.3(97.3) &0.067 ~\\
   0.05        &                 &33.1(31.9)~~~&     0.036~~ &97.8(95.8)&0.020~ \\
  \hline
   Z        &                 &69.5~~~& &96.9~~~~~&\\
    \hline\hline
   0.15    &      $\chi^{2}(2)$                  &46.0(43.2)~~   &0.061 ~~&105.0(97.6) &0.071 ~\\
   0.05        &                 &29.0(28.2)~~~&     0.028~~ &98.1(95.6) &0.026~ \\
  \hline
   Z        &                 &69.7~~~& &96.5~~~~~&\\
\hline\hline
   0.15    &      $t(5)$                 &43.4(39.7)~~   &0.085 ~~&110.2(98.2) &0.109 ~\\
   0.05        &                 &20.2(19.6)~~~&     0.030~~ &99.9(96.6)&0.033~ \\
  \hline
   Z        &                 &67.6~~~& &96.6~~~~~&\\
\hline
 \end{tabular}
\end{center}
}
\end{table}


\section{Proofs}
\setcounter{equation}{0}

 Throughout, we let $C$ denote positive constant, and its value may be different in different contexts. For two
real sequences $\{a_{n}\}$ and $\{b_{n}\}$, write $a_{n} = O(b_{n})$ if there exists a constant $C$ such that $|a_{n}| \leq C|b_{n}|$ holds for
large $n$, $a_{n} = o(b_{n})$ if $\lim_{n\rightarrow\infty}a_{n}/b_{n} = 0$. We denote by $|\cdot|$ the $d$-dimensional Euclidean norm in
$\textbf{R}^{d}$, $d\geq 1$.\\\vspace{-0.2cm}

\noindent{\bf Proof of Theorem 2.1.} (i) ({\bf Sufficiency}) By $\sum_{k=1}^{n}e^{ik\omega_{j}}=0$ for $1\leq j\leq
 q$, we can assume that $\ep X_{1}=0$. Also, for convenience, we assume $\sigma^{2}=1$. For $y\in[-\log q, ~c\log n]$,  set $x=\sqrt{y+\log q}$.
 We start with trunction of $X_{k}$ at two levels. Let $\varepsilon_{n}=(\log n)^{-1}$ and $\varepsilon>0$
 be a small number which will be specified later. Define
\begin{eqnarray}\label{p9}
&&\widetilde{X}_{k}=X^{'}_{k}-\ep X^{'}_{k},\quad X^{'}_{k}=X_{k}I\{|X_{k}|\leq \varepsilon\sqrt{n}x\}, \cr && \widehat{X}_{k}=X^{''}_{k}-\ep
X^{''}_{k},\quad X^{''}_{k}=X_{k}I\{|X_{k}|\leq \varepsilon_{n}\sqrt{n}/x\},\cr
&&I^{'}_{n}(\omega_{j})=\frac{1}{n}\Big{|}\sum_{k=1}^{n}X^{'}_{k}e^{ik\omega_{j}}\Big{|}^{2}=
\frac{1}{n}\Big{|}\sum_{k=1}^{n}\widetilde{X}_{k}e^{ik\omega_{j}}\Big{|}^{2},\cr
&&I^{''}_{n}(\omega_{j})=\frac{1}{n}\Big{|}\sum_{k=1}^{n}X^{''}_{k}e^{ik\omega_{j}}\Big{|}^{2}=
\frac{1}{n}\Big{|}\sum_{k=1}^{n}\widehat{X}_{k}e^{ik\omega_{j}}\Big{|}^{2}.
\end{eqnarray}
 Then we have
\begin{eqnarray}\label{p1}
&&\Big{|}\pr\Big{(}\max_{1\leq j\leq q}I_{n}(\omega_{j})\geq x^{2}\Big{)}- \pr\Big{(}\max_{1\leq j\leq q}I^{'}_{n}(\omega_{j})\geq
x^{2}\Big{)}\Big{|}\cr &&\quad\leq n\pr\Big{(}|X_{1}|\geq \varepsilon\sqrt{n}x\Big{)}.
\end{eqnarray}
 The independence between $X_{l}$ and $\{X_{k}, k\neq l\}$, $1\leq l\leq n$,  implies that
\begin{eqnarray}\label{p2}
&&\pr\Big{(}\max_{1\leq j\leq q}I^{'}_{n}(\omega_{j})\geq x^{2}\Big{)}-\pr\Big{(}\max_{1\leq j\leq q}I^{''}_{n}(\omega_{j})\geq x^{2}\Big{)}\cr
&&\leq \pr\Big{(}\max_{1\leq j\leq q}\Big{|}\sum_{k=1}^{n}X^{'}_{k}e^{ik\omega_{j}}\Big{|}\geq \sqrt{n}x, \cup_{l=1}^{n}\{|X_{l}|>
\varepsilon_{n}\sqrt{n}/x\}\Big{)}\cr &&\leq\sum_{l=1}^{n}\pr\Big{(}\max_{1\leq j\leq q}\Big{|}\sum_{k=1,k\neq
l}^{n}X^{'}_{k}e^{ik\omega_{j}}\Big{|}\geq (1-\varepsilon)\sqrt{n}x\Big{)}\pr\Big{(} |X_{l}|> \varepsilon_{n}\sqrt{n}/x\Big{)}\cr && \leq
n\pr\Big{(}\max_{1\leq j\leq q}\Big{|}\sum_{k=1}^{n}X^{'}_{k}e^{ik\omega_{j}}\Big{|}\geq (1-2\varepsilon)\sqrt{n}x\Big{)} \pr\Big{(} |X_{1}|>
\varepsilon_{n}\sqrt{n}/x\Big{)}\cr && \leq n\pr\Big{(}\max_{1\leq j\leq q}\Big{|}\sum_{k=1}^{n}X^{''}_{k}e^{ik\omega_{j}}\Big{|}\geq
(1-2\varepsilon)\sqrt{n}x\Big{)} \pr\Big{(} |X_{1}|> \varepsilon_{n}\sqrt{n}/x\Big{)}\cr &&\quad+\Big{(}n \pr\Big{(} |X_{1}|>
\varepsilon_{n}\sqrt{n}/x\Big{)}\Big{)}^{2}\cr
 &&=:A+\Big{(}n \pr\Big{(} |X_{1}|>
\varepsilon_{n}\sqrt{n}/x\Big{)}\Big{)}^{2}.
\end{eqnarray}
To estimate $A$, we need Lemma 4.2 of Lin and Liu (2009a). The proof is given in Lin and Liu (2009b), pages 23-25 and the constants $n_{0}$, $C$,
$c_{1,1},\cdots,c_{1,3}$ below are specified in  Lin and Liu (2009b).

\begin{lemma}\label{le5-2}[Lin and Liu (2009a,b)]  Let $\xi_{n,1},\cdots, \xi_{n,k_{n}}$ be independent random
vectors with mean zero
 and values in $\textbf{R}^{2d}$, and $S_{n}=\sum_{i=1}^{k_{n}}X_{n,i}$. Assume that $|\xi_{n,k}|\leq
 c_{n}B^{1/2}_{n}$, $1\leq k\leq k_{n}$, for some $c_{n}\rightarrow 0$,
 $B_{n}\rightarrow\infty$ and
 \begin{equation*}
 \Big{\|}B^{-1}_{n}\Cov(\xi_{n,1}+\cdots+\xi_{n,k_{n}})-I_{2d}\Big{\|}\leq C_{0}c^{2}_{n},
 \end{equation*}
 where $I_{2d}$ is a $2d\times 2d$ identity matrix and $C_{0}$ is a positive constant. Suppose that
$
 \beta_{n}:=B^{-3/2}_{n}\sum_{k=1}^{k_{n}}\ep|\xi_{n,k}|^{3}\rightarrow
 0.
 $
 Then for all $n\geq n_{0}$ ($n_{0}$ is given below)
 \begin{align*}
& |\pr(\|S_{n}\|_{d}\geq x)-\pr(\|N\|_{d}\geq
x/B^{1/2}_{n})|\\
&\leq o(1)\pr(\|N\|_{d}\geq x/B^{1/2}_{n})+C\Big{(}\exp\Big{(}-\frac{\delta^{2}_{n}\min(c^{-2}_{n}, \beta^{-2/3}_{n})}{16d}\Big{)}
+\exp\Big{(}\frac{Cc^{2}_{n}}{\beta^{2}_{n}\log \beta_{n}}\Big{)}\Big{)},
 \end{align*}
uniformly for $x\in[B^{1/2}_{n}, \delta_{n}\min(c^{-1}_{n}, \beta^{-1/3}_{n})B^{1/2}_{n}]$, with any $\delta_{n}\rightarrow 0$ and
$\delta_{n}\min(c^{-1}_{n}, \beta^{-1/3}_{n})\rightarrow \infty$. $N$ is a centered normal random vector with covariance matrix $I_{2d}$.
$\|\cdot\|_{d}$ is defined by $\|z\|_{d}=\min\{(x^{2}_{i}+y^{2}_{i})^{1/2}: 1\leq i\leq
  d\}$, $z=(x_{1}, y_{1},\cdots, x_{d},y_{d})$. $\|\cdot\|$ is the spectral norm for the matrix.
  $C$ is a positive constant which  depends only on $d$ and $C_{0}$. $o(1)$ is bounded by $c_{1,1}(\delta^{3}_{n}+\beta_{n}+c_{n})$, $c_{1,1}$ is
  a positive constant  depending only on $d$.
  \begin{eqnarray*}
n_{0}=\min\Big{\{}n: \forall k\geq n,~c^{2}_{k}\leq \frac{\min(C^{-1}_{0},8^{-1})}{2},~\delta_{k}\leq c_{1,2}\min(C^{-2}_{0},1), ~\beta_{k}\leq
  c_{1,3}\Big{\}},
  \end{eqnarray*}
  where $c_{1,2}$ and $c_{1,3}$ are positive constants depending only on $d$.
\end{lemma}

 Let $d\geq 1$ be a fixed integer and
\begin{eqnarray*}
&&\textbf{Y}_{k}:=\textbf{Y}_{k}(\omega_{i_{1}},\cdots,\omega_{i_{d}})\cr
&&\quad\quad=\widehat{X}_{k}\Big{(}\cos(k\omega_{i_{1}}),\sin(k\omega_{i_{1}}),\cdots, \cos(k\omega_{i_{d}}),\sin(k\omega_{i_{d}})\Big{)},\cr
&&1\leq k\leq n,\quad 1\leq i_{1}<\cdots<i_{d}\leq q.
\end{eqnarray*}
Let $\|z\|_{d}=\min\{(x^{2}_{k}+y^{2}_{k})^{1/2}:1\leq k\leq d\}$ with $z=(x_{1},y_{1},\cdots,x_{d},y_{d})$. By  the facts that for $1\leq j,l\leq
q$,
\begin{eqnarray}\label{eq:fat}
&&\sum_{k=1}^{n}\cos^{2}(\omega_{j}k)=n/2,~\sum_{k=1}^{n}\sin^{2}(\omega_{j}k)=n/2,~\sum_{k=1}^{n}\cos(\omega_{j}k)\sin(\omega_{l}k)=0,
\end{eqnarray}
we can see that
\begin{eqnarray}\label{p10}
\Big{\|}n^{-1}\Cov\Big{(}\sum_{k=1}^{n}\textbf{Y}_{k}\Big{)}-\frac{1}{2} \textbf{I}_{2d}\Big{\|}\leq \ep X^{2}_{1}I\{|X_{1}|\geq
\varepsilon_{n}\sqrt{n}/x\},
\end{eqnarray}
where $\textbf{I}_{2d}$ is $2d\times 2d$ identity matrix and $\|\cdot\|$ is the spectral norm. It is easy to see that (\ref{th1}) implies $\ep
|X_{1}|^{p}<\infty$ for any $2<p<2c+2$. Thus we have
\begin{eqnarray}\label{p11}
&&|\textbf{Y}_{k}|\leq 2d\varepsilon_{n}\sqrt{n}/x\mbox{~and~} n^{-3/2}\sum_{k=1}^{n}\ep|\textbf{Y}_{k}|^{3}\leq
C\max(n^{1-p/2}x^{-3+p},n^{-1/2}).
\end{eqnarray}
Letting $c_{n}=2\varepsilon_{n}(\log q)^{-1/2}$, $B_{n}=n$ in Lemma \ref{le5-2} and $\delta_{n}\log n\rightarrow\infty$ with $\delta_{n}$ being
defined in Lemma \ref{le5-2}, we have for any $0\leq \eta<1$,
\begin{eqnarray}\label{p3}
\frac{\pr\Big{(}\Big{\|}\sum_{k=1}^{n}\textbf{Y}_{k}\Big{\|}_{d}\geq (1-\eta)\sqrt{n}x\Big{)}}{
q^{-d(1-\eta)^{2}}\exp(-d(1-\eta)^{2}y)}\rightarrow 1,
\end{eqnarray}
uniformly in $x\in [\sqrt{\log q}, \sqrt{(c+1)\log n}]$. Observing that
$$I^{''}_{n}(\omega_{j})=\frac{1}{n}\Big{(}\Big{\|}\sum_{k=1}^{n}\textbf{Y}_{k}(\omega_{j})\Big{\|}_{1}\Big{)}^{2},$$ by (\ref{p3}),
\begin{eqnarray*}
\pr\Big{(}\max_{1\leq j\leq q}I^{''}_{n}(\omega_{j})\geq x^{2}\Big{)}\leq\sum_{j=1}^{q}\pr\Big{(}
\Big{\|}\sum_{k=1}^{n}\textbf{Y}_{k}(\omega_{j})\Big{\|}_{1}\geq \sqrt{n}x\Big{)}\leq (1+o(1))\exp(-y)
\end{eqnarray*}
uniformly in $y\in [0, c\log n]$.
  So, combining (\ref{p1})-(\ref{p3}), we prove that
\begin{eqnarray}\label{p12}
\pr\Big{(}\max_{1\leq j\leq q}I_{n}(\omega_{j})\geq x^{2}\Big{)}&\leq& (1+o(1))\exp(-y)\cr & &+Cn^{1+4\varepsilon} \exp(-(1-2\varepsilon)^{2}y)
 \pr\Big{(} |X_{1}|> \varepsilon_{n}\sqrt{n}/x\Big{)}\cr
 & &+\Big{(}n \pr\Big{(} |X_{1}|>
\varepsilon_{n}\sqrt{n}/x\Big{)}\Big{)}^{2}+n\pr\Big{(}|X_{1}|\geq \varepsilon\sqrt{n}x\Big{)},
\end{eqnarray}
uniformly in $x\in [\sqrt{\log q}, \sqrt{(c+1)\log n}]$. To establish the lower bound, we observe that
\begin{eqnarray}\label{p4}
\pr\Big{(}\max_{1\leq j\leq q}I_{n}(\omega_{j})\geq x^{2}\Big{)} &\geq& \pr\Big{(}\max_{1\leq j\leq q}I^{''}_{n}(\omega_{j})\geq x^{2},
\cap_{k=1}^{n}\{|X_{k}|\leq \varepsilon_{n}\sqrt{n}/x\}\Big{)}\cr &=& \pr\Big{(}\max_{1\leq j\leq q}I^{''}_{n}(\omega_{j})\geq x^{2}\Big{)}\cr &
&-\pr\Big{(}\max_{1\leq j\leq q}I^{''}_{n}(\omega_{j})\geq x^{2}, \cup_{k=1}^{n}\{|X_{k}|\geq \varepsilon_{n}\sqrt{n}/x\}\Big{)}.
\end{eqnarray}
Similarly to (\ref{p2}) and by  (\ref{p3}), we have
\begin{eqnarray}\label{p5}
&&\pr\Big{(}\max_{1\leq j\leq q}I^{''}_{n}(\omega_{j})\geq x^{2}, \cup_{k=1}^{n}\{|X_{k}|\geq \varepsilon_{n}\sqrt{n}/x\}\Big{)}\cr &&\quad \leq
n\pr\Big{(}\max_{1\leq j\leq q}\Big{|}\sum_{k=1}^{n}X^{''}_{k}e^{ik\omega_{j}}\Big{|}\geq(1-2\varepsilon)\sqrt{n}x\Big{)}\pr(|X_{1}|\geq
\varepsilon_{n}\sqrt{n}/x)\cr
 &&\quad \leq
Cn^{1+4\varepsilon} \exp(-(1-2\varepsilon)^{2}y)
 \pr\Big{(} |X_{1}|> \varepsilon_{n}\sqrt{n}/x\Big{)}
\end{eqnarray}
uniformly in $x\in [\sqrt{\log q}, \sqrt{(c+1)\log n}]$. For the first term in (\ref{p4}), we have
\begin{eqnarray}\label{p6}
\pr\Big{(}\max_{1\leq j\leq q}I^{''}_{n}(\omega_{j})\geq x^{2}\Big{)}\geq\sum_{j=1}^{q}\pr(A_{j}) - \sum_{1\leq i<j\leq q}\pr(A_{i}A_{j}),
\end{eqnarray}
where $ A_{j}=\{I^{''}_{n}(\omega_{j})\geq x^{2}\}.$  Applying $d=1,2$ in (\ref{p3}) respectively,  we obtain
\begin{eqnarray*}
\pr(A_{i})=(1+o(1))q^{-1}\exp(-y),\quad\pr(A_{i}A_{j})\leq Cn^{-2}\exp(-2y),
\end{eqnarray*}
uniformly in $1\leq i,j\leq q$ and $y\in [0, c\log n]$. These two inequalities together with (\ref{p4})-(\ref{p6}) yield that
\begin{eqnarray}\label{p13}
\pr\Big{(}\max_{1\leq j\leq q}I_{n}(\omega_{j})\geq x^{2}\Big{)}&\geq& (1+o(1))\exp(-y)-C\exp(-2y)\cr & & -Cn^{1+4\varepsilon}
\exp(-(1-2\varepsilon)^{2}y)
 \pr\Big{(} |X_{1}|> \varepsilon_{n}\sqrt{n}/x\Big{)},
\end{eqnarray}
uniformly in $y\in [0, c\log n]$. It is easy to see that (\ref{th1}) implies that
\begin{eqnarray*}
\pr\Big{(} |X_{1}|> \varepsilon_{n}\sqrt{n}/x\Big{)}\leq C\frac{(\log n)^{4c+4}}{n^{c+1}}
\end{eqnarray*}
for $\varepsilon_{n}=(\log n)^{-1}$ and $x\in [\sqrt{\log q}, \sqrt{(c+1)\log n}]$. Hence, by (\ref{p12}), (\ref{p13}) and for $\varepsilon$
sufficiently small, we have for any $M>0$,
\begin{eqnarray}\label{eq:g3}
&&\limsup_{n\rightarrow\infty}\sup_{M\leq y\leq c\log n}\Big{|}\frac{\pr\Big{(}\max_{1\leq j\leq q}I_{n}(\omega_{j})-\log q\geq
y\Big{)}}{1-\exp(-\exp(-y))}-1\Big{|}\leq Ce^{-M}.
\end{eqnarray}
By (\ref{eq:f}), we have for any fixed $y\in \textbf{R}$,
\begin{eqnarray*}
\pr\Big{(}\max_{1\leq j\leq q}I_{n}(\omega_{j})-\log q\geq y\Big{)}\rightarrow 1-e^{-e^{-y}}.
\end{eqnarray*}
Thus, it follows that
\begin{eqnarray}\label{eq:g4}
&&\limsup_{n\rightarrow\infty}\sup_{-\log q\leq y\leq M}\Big{|}\frac{\pr\Big{(}\max_{1\leq j\leq q}I_{n}(\omega_{j})-\log q\geq
y\Big{)}}{1-\exp(-\exp(-y))}-1\Big{|}=0.
\end{eqnarray}
This proves (i) by (\ref{eq:g3}) and (\ref{eq:g4}).

(ii) ({\bf Necessity}) By  (\ref{th2}), we have
\begin{eqnarray*}
\pr\Big{(}\max_{1\leq j\leq q}\Big{|}\sum_{k=1}^{n}X_{k}e^{ik\omega_{j}}\Big{|}\geq \sigma\sqrt{(1+c)n\log q}\Big{)}\leq Cn^{-c}.
\end{eqnarray*}
This implies that
\begin{eqnarray}\label{p7}
\pr\Big{(}\max_{1\leq j\leq q}\Big{|}\sum_{k=1}^{n}X^{s}_{k}e^{ik\omega_{j}}\Big{|}\geq 2^{-1}\sigma\sqrt{(1+c)n\log q}\Big{)}\leq Cn^{-c},
\end{eqnarray}
where $X^{s}_{n}=X_{n}-X^{c}_{n}$ and $\{X^{c}_{n}\}$ is an independent copy of $\{X_{n}\}$.  For $z=(z_{1,1},z_{1,2}\ldots,$ $z_{q,1},
z_{q,2})\in\textbf{R}^{2q}$, let $$\|z\|_{\max}=\max_{1\leq j\leq q}\sqrt{z^{2}_{j,1}+z^{2}_{j,2}}.$$ For $1\leq k\leq n$, let
$$\texttt{D}_{k}=(X^{s}_{k}\cos(k \omega_{1}), X^{s}_{k}\sin(k \omega_{1}),\cdots, X^{s}_{k}\cos(k \omega_{q}), X^{s}_{k}\sin(k \omega_{q})).$$
Then it is easy to see that $\max_{1\leq j\leq q}\Big{|}\sum_{k=1}^{n}X^{s}_{k}e^{ik\omega_{j}}\Big{|}=\|\sum_{k=1}^{n}\texttt{D}_{k}\|_{\max}$.
 By L\'{e}vy's inequality in a Banach space (cf. Ledoux and Talagrand (1991), p. 47) and (\ref{p7}), we have
\begin{eqnarray}\label{p8}
\pr\Big{(}\max_{1\leq k\leq n}\|\texttt{D}_{k}\|_{\max}\geq 2^{-1}\sigma\sqrt{(1+c)n\log q}\Big{)}\leq Cn^{-c}.
\end{eqnarray}
Observing that $\|\texttt{D}_{k}\|_{\max}=|X^{s}_{k}|$,  we have by (\ref{p8}),
\begin{eqnarray*}
1-\Big{(}1-\pr\Big{(}|X^{s}_{1}|\geq 2^{-1}\sigma\sqrt{(1+c)n\log q}\Big{)}\Big{)}^{n}\leq Cn^{-c},
\end{eqnarray*}
which implies that
\begin{eqnarray}\label{eq:g5}
1-\exp\Big{(}-n\pr\Big{(}|X^{s}_{1}|\geq 2^{-1}\sigma\sqrt{(1+c)n\log q}\Big{)}\Big{)}\leq Cn^{-c}.
\end{eqnarray}
Now (\ref{th3}) follows from (\ref{eq:g5}) and some elementary  calculations. $\Box$\vspace{0.5cm}

\noindent{\bf Proof of Theorem 2.2.} For any  $\varepsilon^{'}_{n}\rightarrow 0$, let $x=\sqrt{y+\log q}\in [4, \varepsilon^{'}_{n}n^{1/6}]$ and
$\varepsilon_{n}=\varepsilon^{'1/4}_{n}$. Define
\begin{eqnarray*}
X^{'}_{k}=X_{k}I\{|X_{k}|\leq \varepsilon_{n}\sqrt{n}/x\},\quad
I^{'}_{n}(\omega_{j})=\frac{1}{n}\Big{|}\sum_{k=1}^{n}X^{'}_{k}e^{ik\omega_{j}}\Big{|}^{2}.
\end{eqnarray*}
Then, by Lemma 4.2 below (taking $\textbf{N}_{l}=\emptyset$ and $d=1$), we have
\begin{eqnarray*}
\pr\Big{(}\max_{1\leq j\leq q}I^{'}_{n}(\omega_{j})\geq x^{2}\Big{)}\leq
\sum_{j=1}^{q}\pr\Big{(}\Big{|}\sum_{k=1}^{n}X^{'}_{k}e^{ik\omega_{j}}\Big{|}\geq \sqrt{n}x\Big{)}\leq (1+o(1))e^{-y}
\end{eqnarray*}
uniformly in $x\in[4, \varepsilon^{'}_{n}n^{1/6}]$. Since $\ep e^{t_{0}|X_{1}|}<\infty$, we have $n\pr\Big{(} |X_{1}|>
\varepsilon_{n}\sqrt{n}/x\Big{)}=o(1)e^{-x^{2}}.$ Therefore,
\begin{eqnarray*}
\pr\Big{(}\max_{1\leq j\leq q}I_{n}(\omega_{j})\geq x^{2}\Big{)}&\leq& \pr\Big{(}\max_{1\leq j\leq q}I^{'}_{n}(\omega_{j})\geq x^{2}\Big{)}+
n\pr\Big{(} |X_{1}|> \varepsilon_{n}\sqrt{n}/x\Big{)}\cr &\leq& (1+o(1))e^{-y}
\end{eqnarray*}
uniformly in $x\in[4, \varepsilon^{'}_{n}n^{1/6}]$. Similarly, we have
\begin{eqnarray*}
\pr\Big{(}\max_{1\leq j\leq q}I_{n}(\omega_{j})\geq x^{2}\Big{)}\geq (1+o(1))e^{-y}-Ce^{-2y}
\end{eqnarray*}
uniformly in $x\in[4, \varepsilon^{'}_{n}n^{1/6}]$. The remaining proof follows similar arguments as in the proof of Theorem 2.1.
$\Box$\vspace{0.5cm}

\noindent{\bf Proof of Theorem 2.3.} Recall that $\sum_{k=1}^{n}e^{ik\omega_{j}}=0$ for $1\leq j\leq
 q$. Without loss of generality, we can assume that $\ep
 X_{1}=0$ and $\ep X^{2}_{1}=1$.
By the fact $e^{ik\omega_{j_{1}}}=e^{-ik\omega_{j_{2}}}$ for $j_{1}+j_{2}=n$, we have
\begin{eqnarray*}
\Big{|}\sum_{k=1}^{n}X_{k}e^{ik\omega_{j_{1}}}\Big{|}^{2} =\Big{|}\sum_{k=1}^{n}X_{k}e^{ik\omega_{j_{2}}}\Big{|}^{2}
\end{eqnarray*}
for $j_{1}+j_{2}=n$. Hence, when $n$ is odd,
\begin{eqnarray}\label{f1}
q^{-1}\sum_{j=1}^{q}I_{n}(\omega_{j})&=&(2q)^{-1}\sum_{j=1}^{q}\Big{(}I_{n}(\omega_{j})+ I_{n}(\omega_{n-j})\Big{)}\cr
&=&(2q)^{-1}\sum_{j=1}^{n}I_{n}(\omega_{j})-\frac{n}{n-1}(\overline{X})^{2},
\end{eqnarray}
where $\overline{X}=n^{-1}\sum_{k=1}^{n}X_{k}$. Moreover,
\begin{eqnarray*}
\sum_{j=1}^{n}I_{n}(\omega_{j})=\sum_{k=1}^{n}X^{2}_{k}+2n^{-1}\sum_{k=2}^{n}X_{k}\sum_{i=1}^{k-1}X_{i} \sum_{j=1}^{n}w_{k,i,j},
\end{eqnarray*}
where $w_{k,i,j}=\cos(\omega_{j}k)\cos(\omega_{j}i)+\sin(\omega_{j}k)\sin(\omega_{j}i)$. Note that $\omega_{j}k=\omega_{k}j$ and
$\omega_{j}i=\omega_{i}j$. Since $|\sum_{l=1}^{n}e^{il\lambda}|=|\sin(\lambda n/2)|/|\sin(\lambda/2)|$ when $\lambda/\pi$ is not an integer, we
 get $\sum_{j=1}^{n}w_{k.i.j}=\sum_{j=1}^{n}\cos((\omega_{k}-\omega_{i})j)=0$. So, (\ref{f1}) implies that, when $n$ is odd,
\begin{eqnarray*}
q^{-1}\sum_{j=1}^{q}I_{n}(\omega_{j})=(n-1)^{-1}\Big{(}\sum_{k=1}^{n}X^{2}_{k}-n(\overline{X})^{2}\Big{)}.
\end{eqnarray*}
Similarly, when $n$ is even, we have
\begin{eqnarray}\label{f2}
q^{-1}\sum_{j=1}^{q}I_{n}(\omega_{j})=(n-2)^{-1}\Big{(} \sum_{k=1}^{n}X^{2}_{k}-n(\overline{X})^{2}-n(\overline{X}^{'})^{2}\Big{)},
\end{eqnarray}
where $\overline{X}^{'}=n^{-1}\sum_{k=1}^{n}(-1)^{k}X_{k}$. By the self-normalized moderate deviation Theorem 3.1 in Shao (1997), we have
\begin{eqnarray*}
\pr\Big{(}\Big{|}\sum_{k=1}^{n}(-1)^{k}X_{k}\Big{|}\geq n^{1/3}V_{n}\Big{)}\leq Ce^{-n^{2/3}/4},\quad
\pr\Big{(}\Big{|}\sum_{k=1}^{n}X_{k}\Big{|}\geq n^{1/3}V_{n}\Big{)}\leq Ce^{-n^{2/3}/4},
\end{eqnarray*}
where $V^{2}_{n}=\sum_{k=1}^{n}X^{2}_{k}$. In view of (\ref{f1}) and (\ref{f2}), it suffices to show that
\begin{eqnarray}\label{eq:g6}
\lim_{n\rightarrow\infty}\frac{\pr\Big{(}\frac{\max_{1\leq j\leq
q}\Big{|}\sum_{k=1}^{n}X_{k}e^{ik\omega_{j}}\Big{|}^{2}}{\sum_{k=1}^{n}X^{2}_{k}}-\log q\geq y\Big{)}}{1-e^{-e^{-y}}}=1
\end{eqnarray}
uniformly in $y\in [-\log q, ~o(n^{1/3})).$

  Let $\lambda=\lambda_{n}$ be a positive number which will be specified
  later. Let $\textbf{H}$ be a subset of $\{1,\cdots,n\}$. Put
\begin{eqnarray*}
&& X^{'}_{k}=X_{k}I\{|X_{k}|\leq \lambda\},\quad 1\leq k\leq n,\\
 &&M_{n}=\max_{1\leq j\leq
q}\Big{|}\sum_{k=1}^{n}X_{k}e^{ik\omega_{j}}\Big{|},\quad \widetilde{M}_{n}=\max_{1\leq j\leq
q}\Big{|}\sum_{k=1}^{n}X^{'}_{k}e^{ik\omega_{j}}\Big{|},\\
&&M^{(\textbf{H})}_{n}=\max_{1\leq j\leq q}\Big{|}\sum_{k=1,k\notin \textbf{H}}^{n}X_{k}e^{ik\omega_{j}}\Big{|},\quad
\widetilde{M}^{(\textbf{H})}_{n}=\max_{1\leq j\leq q}\Big{|}\sum_{k=1,k\notin
\textbf{H}}^{n}X^{'}_{k}e^{ik\omega_{j}}\Big{|},\\
&&\widetilde{V}_{n}=(\sum_{k=1}^{n}X^{'2}_{k})^{1/2},\quad V^{(\textbf{H})}_{n}=(\sum_{k=1,k\notin \textbf{H}}^{n}X^{2}_{k})^{1/2},\quad
\widetilde{V}^{(\textbf{H})}_{n}=(\sum_{k=1,k\notin\textbf{H}}^{n}X^{'2}_{k})^{1/2}.
\end{eqnarray*}
Noting that for any real numbers $s$ and $t$ and nonnegative number $c$ and $x\geq 1$,
\begin{eqnarray}\label{gin}
\{s+t\geq x\sqrt{c+t^{2}}\}\subset\{s\geq\sqrt{x^{2}-1}\sqrt{c}\},
\end{eqnarray}
(see p.2181 in Jing, Shao and Wang (2003)), we have
\begin{eqnarray}\label{eq:g9-1}
\pr\Big{(}M_{n}\geq xV_{n}\Big{)}&\leq& \pr\Big{(}\widetilde{M}_{n}\geq x\widetilde{V}_{n}\Big{)}+\sum_{j=1}^{n}\pr\Big{(}M_{n}\geq
xV_{n},X_{j}\neq X^{'}_{j}\Big{)}\cr &\leq&\pr\Big{(}\widetilde{M}_{n}\geq
x\widetilde{V}_{n}\Big{)}+\sum_{j=1}^{n}\pr\Big{(}M^{(j)}_{n}+|X_{j}|\geq xV_{n},X_{j}\neq X^{'}_{j}\Big{)}\cr
&\leq&\pr\Big{(}\widetilde{M}_{n}\geq x\widetilde{V}_{n}\Big{)}+\sum_{j=1}^{n}\pr\Big{(}M^{(j)}_{n}\geq \sqrt{x^{2}-1}V^{(j)}_{n},X_{j}\neq
X^{'}_{j}\Big{)}\cr &=&\pr\Big{(}\widetilde{M}_{n}\geq x\widetilde{V}_{n}\Big{)}+\sum_{j=1}^{n}\pr\Big{(}|X_{j}|\geq
\lambda\Big{)}\pr\Big{(}M^{(j)}_{n}\geq \sqrt{x^{2}-1}V^{(j)}_{n}\Big{)}.\quad\quad
\end{eqnarray}
Repeating the above arguments $m$ times with $m=[x^{2}/2]$, we have for $x>4$,
\begin{eqnarray}\label{eq:g9-3}
&&\sum_{j_{1}=1}^{n}\pr\Big{(}|X_{j_{1}}|\geq \lambda\Big{)}\pr\Big{(}M^{(j_{1})}_{n}\geq \sqrt{x^{2}-1}V^{(j_{1})}_{n}\Big{)}\cr
&&\quad\leq\sum_{j_{1}=1}^{n}\pr\Big{(}|X_{j_{1}}|\geq \lambda\Big{)}\pr\Big{(}\widetilde{M}^{(j_{1})}_{n}\geq
\sqrt{x^{2}-1}\widetilde{V}^{(j_{1})}_{n}\Big{)}\cr &&\quad\quad+\sum_{j_{1}=1}^{n}\sum_{j_{2}=1}^{n}\pr\Big{(}|X_{j_{1}}|\geq
\lambda\Big{)}\pr\Big{(}|X_{j_{2}}|\geq \lambda\Big{)}\pr\Big{(}M^{(j_{1},j_{2})}_{n}\geq \sqrt{x^{2}-2}V^{(j_{1},j_{2})}_{n}\Big{)}\cr
&&\quad\leq \sum_{l=1}^{m}\widetilde{Z}_{l}+Z_{m+1},
\end{eqnarray}
where
\begin{eqnarray*}
&&\widetilde{Z}_{l}=\sum_{j_{1}=1}^{n}\cdots\sum_{j_{l}=1}^{n}\Big{[}\prod_{k=1}^{l}\pr\Big{(}|X_{j_{k}}|\geq \lambda\Big{)}\Big{]}\times
\pr\Big{(}\widetilde{M}^{(j_{1},\cdots,j_{l})}_{n}\geq \sqrt{x^{2}-l}\widetilde{V}^{(j_{1},\cdots,j_{l})}_{n}\Big{)},\cr
&&Z_{m+1}=\sum_{j_{1}=1}^{n}\cdots\sum_{j_{m+1}=1}^{n}\Big{\{}\Big{[}\prod_{k=1}^{m+1}\pr\Big{(}|X_{j_{k}}|\geq \lambda\Big{)}\Big{]}\cr
&&\qquad\qquad\qquad\qquad\quad\quad~\times \pr\Big{(}M^{(j_{1},\cdots,j_{m+1})}_{n}\geq
\sqrt{x^{2}-m-1}V^{(j_{1},\cdots,j_{m+1})}_{n}\Big{)}\Big{\}}.
\end{eqnarray*}
For $\varepsilon^{'}_{n}\rightarrow 0$ and $4\leq x\leq \varepsilon^{'}_{n}n^{1/6}$, let  $\lambda=\varepsilon_{n}\sqrt{n}/x$, where
$\varepsilon_{n}=\varepsilon^{'1/4}_{n}$. Then
\begin{eqnarray}\label{eq:g9-2}
Z_{m+1}\leq \Big{(}n\pr\Big{(}|X_{1}|\geq \lambda\Big{)}\Big{)}^{m+1} \leq e^{-m\log q_{n}}=o(1)e^{-x^{2}},
\end{eqnarray}
where
\begin{eqnarray*}
q_{n}=\Big{(}\varepsilon^{'3}_{n}\varepsilon^{-3}_{n}\ep |X_{1}|^{3}I\{|X_{1}|\geq\lambda\}\Big{)}^{-1}\rightarrow\infty
\end{eqnarray*}
as $n\rightarrow\infty$. Let $0\leq l\leq m=[x^{2}/2]$ and  $\textbf{N}_{l}=\{j_{1},\cdots,j_{l}\}\subset \{1,\cdots,n\}$.  Define
\begin{eqnarray*}
&&\textbf{Y}^{'}_{k}=X^{'}_{k}\Big{(}\cos(k\omega_{i_{1}}),\sin(k\omega_{i_{1}}),\cdots, \cos(k\omega_{i_{d}}),\sin(k\omega_{i_{d}})\Big{)},\cr
&&1\leq k\leq n,\quad d\geq 1,\quad 1\leq i_{1}<\cdots<i_{d}\leq q.\cr &&
 S^{\textbf{N}_{l}}_{n}=\sum_{k=1,k\notin \textbf{N}_{l}}^{n}\textbf{Y}^{'}_{k},\quad
\overline{S}^{\textbf{N}_{l}}_{n}=\sum_{k=1,k\notin \textbf{N}_{l}}^{n}(\textbf{Y}^{'}_{k}-\ep \textbf{Y}^{'}_{k}).
\end{eqnarray*}
 To estimate $\widetilde{Z}_{l}$, we need the following lemma.
\begin{lemma}\label{le4.2}  Suppose that $\ep |X_{1}|^{3}<\infty$ and $0\leq x\leq
\varepsilon^{'}_{n}n^{1/6}$, where $\varepsilon^{'}_{n}\rightarrow 0$ is any sequence of constants. Let $0<\varepsilon_{n}\rightarrow 0$ and
$\varepsilon_{n}=\varepsilon^{'1/4}_{n}$. (i) If $\lambda=\varepsilon_{n}\sqrt{n}/x$, then we have
\begin{eqnarray}\label{eq8}
\lim_{n\rightarrow\infty}\frac{\pr\Big{(}\Big{\|}S^{\textbf{N}_{l}}_{n}\Big{\|}_{d}\geq x\sqrt{n}\Big{)}}{e^{-dx^{2}}}= 1
\end{eqnarray}
uniformly in $x\in [4, \varepsilon^{'}_{n}n^{1/6}]$, $ 1\leq i_{1}<\cdots<i_{d}\leq q$ and $0\leq l\leq m$. (ii) If
$\lambda=(\varepsilon_{n}\sqrt{n}/x)^{3/4}$ and $\ep X^{4}_{1}<\infty$, then (\ref{eq8}) holds.
\end{lemma}

\begin{proof}
 Recall that $\sum_{k=1}^{n}e^{ik\omega_{j}}=0$ for $1\leq j\leq
 q$.  We have for $n$ large,
 \begin{eqnarray*}
\Big{|}\sum_{k=1,k\notin \textbf{N}_{l}}^{n}\ep \textbf{Y}^{'}_{k}\Big{|}=\Big{|} \sum_{k\in \textbf{N}_{l}}^{n}\ep \textbf{Y}^{'}_{k}\Big{|}\leq
dx^{2}/2 \leq \varepsilon_{n}\sqrt{n}/x.
 \end{eqnarray*}
It follows that, for $x\in [4, \varepsilon^{'}_{n}n^{1/6}]$,
\begin{eqnarray*}
\pr\Big{(}\Big{\|}\overline{S}^{\textbf{N}_{l}}_{n}\Big{\|}_{d}\geq
x\sqrt{n}+\varepsilon_{n}\sqrt{n}/x\Big{)}&\leq&\pr\Big{(}\Big{\|}S^{\textbf{N}_{l}}_{n}\Big{\|}_{d}\geq x\sqrt{n}\Big{)}\cr &\leq&
\pr\Big{(}\Big{\|}\overline{S}^{\textbf{N}_{l}}_{n}\Big{\|}_{d}\geq x\sqrt{n}-\varepsilon_{n}\sqrt{n}/x\Big{)}.
\end{eqnarray*}
Since $\ep|X_{1}|^{3}<\infty$, we have
\begin{eqnarray*}
|\textbf{Y}^{'}_{k}|\leq 2d \varepsilon_{n}\sqrt{n}/x\mbox{~and~} n^{-3/2}\sum_{k=1}^{n}\ep|\textbf{Y}_{k}|^{3}\leq Cn^{-1/2}.
\end{eqnarray*}
Also, simple calculations show that
\begin{eqnarray*}
\Big{\|}\frac{1}{n}\Cov\Big{(}\overline{S}^{\textbf{N}}_{n}\Big{)}-\textbf{I}_{2d}\Big{\|}\leq \ep X^{2}_{1}I\{|X_{1}|\geq
\lambda\}+C_{d}n^{-1}x^{2}\leq Cxn^{-1/2}\varepsilon^{-1}_{n}\leq C\varepsilon^{2}_{n}/x^{2}
\end{eqnarray*}
for $4\leq x\leq \varepsilon^{'}_{n}n^{1/6}$. By taking $c_{n}=2d\varepsilon_{n}/x$, $B_{n}=\sqrt{n}$ and
$\delta_{n}/\varepsilon_{n}\rightarrow\infty$ in Lemma \ref{le5-2}, we have
\begin{eqnarray*}
\frac{\pr\Big{(}\Big{\|}\overline{S}^{\textbf{N}_{l}}_{n}\Big{\|}_{d}\geq x\sqrt{n}\pm\varepsilon_{n}\sqrt{n}/x\Big{)}}{e^{-dx^{2}}}\rightarrow 1,
\end{eqnarray*}
uniformly in $x\in [4, \varepsilon^{'}_{n}n^{1/6}]$. This proves Lemma \ref{le4.2}.
\end{proof}

From Lemma \ref{le4.2}, we have for $0\leq l\leq m=[x^{2}/2]$,
\begin{eqnarray}\label{eq:ft}
&&\pr\Big{(}\widetilde{M}^{(j_{1},\cdots,j_{l})}_{n}\geq \sqrt{x^{2}-l}\widetilde{V}^{(j_{1},\cdots,j_{l})}_{n}\Big{)}\cr
&&\quad\leq\pr\Big{(}\widetilde{M}^{(j_{1},\cdots,j_{l})}_{n}\geq \sqrt{x^{2}-l}\sqrt{n(1-\varepsilon_{n}x^{-2})}\Big{)}\cr
&&\quad\quad+\pr\Big{(}\widetilde{V}^{(j_{1},\cdots,j_{l})}_{n}\leq \sqrt{n(1-\varepsilon_{n}x^{-2})}\Big{)}\cr
&&\quad=(1+o(1))qe^{-x^{2}+l}+\pr\Big{(}\widetilde{V}^{(j_{1},\cdots,j_{l})}_{n}\leq \sqrt{n(1-\varepsilon_{n}x^{-2})}\Big{)},
\end{eqnarray}
uniformly in $x\in [4, \varepsilon^{'}_{n}n^{1/6}]$. By a similar argument as in Hu, Shao and Wang (2009), p. 1193, if $\ep X^{4}_{1}<\infty$, we
have
\begin{eqnarray*}
\pr\Big{(}\widetilde{V}^{(j_{1},\cdots,j_{l})}_{n}\leq
\sqrt{n(1-\varepsilon_{n}x^{-2})}\Big{)}\leq\pr\Big{(}n-m_{l}-\widetilde{V}^{(j_{1},\cdots,j_{l})2}_{n} \geq \varepsilon_{n}nx^{-2}/2\Big{)}\leq
o(1)e^{-x^{2}},
\end{eqnarray*}
where $m_{l}$ is the cardinality of $\{j_{1},\cdots,j_{l}\}$, and hence
\begin{eqnarray}\label{eq:g7}
\widetilde{Z}_{l}\leq Cn\Big{(}n\pr\Big{(}|X_{1}|\geq \lambda\Big{)}\Big{)}^{l}e^{-(x^{2}-l)}\leq Cne^{-x^{2}+l-Cl\log q_{n}}.
\end{eqnarray}
This together with $q_{n}\rightarrow\infty$  shows that
\begin{eqnarray}\label{eq:g9-4}
\sum_{l=1}^{m}\widetilde{Z}_{l}=o(1)ne^{-x^{2}},
\end{eqnarray}
uniformly in $x\in [4, \varepsilon^{'}_{n}n^{1/6}]$. Combining (\ref{eq:g9-1})-(\ref{eq:g9-2}), (\ref{eq:ft}) and (\ref{eq:g9-4}) yields
\begin{eqnarray}\label{p15}
\pr\Big{(}M_{n}\geq xV_{n}\Big{)} \leq (1+o(1))qe^{-x^{2}}
\end{eqnarray}
uniformly in $x\in [4, \varepsilon^{'}_{n}n^{1/6}]$.

We next estimate the lower bound for $\pr\Big{(}M_{n}\geq xV_{n}\Big{)}$. For $\varepsilon^{'}_{n}\rightarrow 0$ and $4\leq x\leq
\varepsilon^{'}_{n}n^{1/6}$, let $\varepsilon_{n}=\max((x/\sqrt{n})^{1/8},\varepsilon^{'1/4}_{n})$ and
$\lambda=(\varepsilon_{n}\sqrt{n}/x)^{3/4}$. Then
\begin{eqnarray}\label{eq:g8}
\pr\Big{(}M_{n}\geq xV_{n}\Big{)}&\geq& \pr\Big{(}\widetilde{M}_{n}\geq x\widetilde{V}_{n}\Big{)}\cr&
&-\sum_{j=1}^{n}\pr\Big{(}\widetilde{M}^{(j)}_{n}\geq \sqrt{x^{2}-1}\widetilde{V}^{(j)}_{n}\Big{)}\pr\Big{(}|X_{j}|\geq \lambda\Big{)}.
\end{eqnarray}
Similarly to (\ref{eq:g7}), we have
\begin{eqnarray}\label{eq9}
\sum_{j=1}^{n}\pr\Big{(}\widetilde{M}^{(j)}_{n}\geq \sqrt{x^{2}-1}\widetilde{V}^{(j)}_{n}\Big{)}\pr\Big{(}|X_{j}|\geq \lambda\Big{)}
=o(1)ne^{-x^{2}}
\end{eqnarray}
uniformly in $x\in [4, \varepsilon^{'}_{n}n^{1/6}]$. For the first term on the right hand side of (\ref{eq:g8}), we have
\begin{eqnarray*}
\pr\Big{(}\widetilde{M}_{n}\geq x\widetilde{V}_{n}\Big{)}&\geq& \pr\Big{(}\widetilde{M}_{n}\geq x\widetilde{V}_{n}, \widetilde{V}^{2}_{n}\leq
n(1+\varepsilon_{n}/x^{2})\Big{)}\cr &\geq& \pr\Big{(}\widetilde{M}_{n}\geq x\sqrt{n(1+\varepsilon_{n}/x^{2})}\Big{)}\cr & &-
\pr\Big{(}\widetilde{M}_{n}\geq\sqrt{n} x, \widetilde{V}^{2}_{n}\geq n(1+\varepsilon_{n}/x^{2})\Big{)}.
\end{eqnarray*}
Define $\textbf{A}=\{\widetilde{V}^{2}_{n}\geq n(1+\varepsilon_{n}/x^{2})\}$. Set
$Y_{k,l}(\theta_{1},\theta_{2})=X^{'}_{k}(\theta_{1}\cos(kw_{l})+\theta_{2}\sin(kw_{l}))$ for any $\theta_{1},\theta_{2}\in \textbf{R}$. Let
\begin{eqnarray*}
&&\Theta_{1}=\{\theta_{1}\geq 0,\theta_{2}\geq 0; \theta^{2}_{1}+\theta^{2}_{2}=1\};~ \Theta_{2}=\{\theta_{1}\geq 0,\theta_{2}<0;
\theta^{2}_{1}+\theta^{2}_{2}=1\};\cr && \Theta_{3}=\{\theta_{1}< 0,\theta_{2}\geq0; \theta^{2}_{1}+\theta^{2}_{2}=1\};~ \Theta_{4}=\{\theta_{1}<
0,\theta_{2}<0; \theta^{2}_{1}+\theta^{2}_{2}=1\}.
\end{eqnarray*}
Then we have
\begin{eqnarray}\label{eq:j10}
&&\pr\Big{(}\widetilde{M}_{n}\geq \sqrt{n}x,\textbf{A}\Big{)}\cr &&\quad\leq\sum_{l=1}^{q}\pr\Big{(} \sup_{(\theta_{1},\theta_{2})\in
\Theta_{1}}\Big{|}\sum_{k=1}^{n} Y_{k,l}(\theta_{1},\theta_{2})\Big{|}\geq \sqrt{n}x,\textbf{A}\Big{)}\cr
 &&\quad\quad+\sum_{l=1}^{q}\pr\Big{(} \sup_{(\theta_{1},\theta_{2})\in
\Theta_{2}}\Big{|}\sum_{k=1}^{n} Y_{k,l}(\theta_{1},\theta_{2})\Big{|}\geq \sqrt{n}x,\textbf{A}\Big{)}\cr &&\quad\quad+\sum_{l=1}^{q}\pr\Big{(}
\sup_{(\theta_{1},\theta_{2})\in \Theta_{3}}\Big{|}\sum_{k=1}^{n} Y_{k,l}(\theta_{1},\theta_{2})\Big{|}\geq \sqrt{n}x,\textbf{A}\Big{)}\cr
&&\quad\quad+\sum_{l=1}^{q}\pr\Big{(} \sup_{(\theta_{1},\theta_{2})\in \Theta_{4}}\Big{|}\sum_{k=1}^{n} Y_{k,l}(\theta_{1},\theta_{2})\Big{|}\geq
\sqrt{n}x,\textbf{A}\Big{)}.
\end{eqnarray}
We only deal with the first term above, while other terms can be proved similarly.
 Let $\theta_{1,i}=i/n^{6}$ for $1\leq i\leq n^{6}$ and
$\theta_{2,i}=\sqrt{1-\theta^{2}_{1,i}}$. We have, for $1\leq l\leq q$,
\begin{eqnarray}\label{eq:j11}
&&\pr\Big{(} \sup_{(\theta_{1},\theta_{2})\in \Theta_{1}}\Big{|}\sum_{k=1}^{n} Y_{k,l}(\theta_{1},\theta_{2})\Big{|}\geq
\sqrt{n}x,\textbf{A}\Big{)}\cr &&\quad\leq \sum_{i=1}^{n^{6}} \pr\Big{(}\Big{|}\sum_{k=1}^{n} Y_{k,l}(\theta_{1,i},\theta_{2,i})\Big{|}\geq
\sqrt{n}x-\sqrt{n}\varepsilon_{n}x^{-1},\textbf{A}\Big{)}\cr & &\quad\quad+ \sum_{i=1}^{n^{6}}
\pr\Big{(}\sup_{(\theta_{1},\theta_{2})\in\Theta_{1},\theta_{1,i-1}\leq \theta_{1}\leq \theta_{1,i}}\Big{|}\sum_{k=1}^{n}
\Big{[}Y_{k,l}(\theta_{1},\theta_{2})\cr &&\qquad\qquad\qquad\qquad\qquad\qquad
\qquad\qquad\qquad-Y_{k,l}(\theta_{1,i},\theta_{2,i})\Big{]}\Big{|}\geq \sqrt{n}\varepsilon_{n}x^{-1}\Big{)}\cr &&\quad
=:\sum_{i=1}^{n^{6}}J_{1,i}+\sum_{i=1}^{n^{6}}J_{2,i}.
\end{eqnarray}
It is easy to see that $\sup_{(\theta_{1},\theta_{2})\in\Theta_{1},\theta_{1,i-1}\leq \theta_{1}\leq
\theta_{1,i}}|Y_{k,l}(\theta_{1},\theta_{2})-Y_{k,l}(\theta_{1,i},\theta_{2,i})|\leq n^{-1}$. Hence $J_{2,i}=0$ for $x\in [4,
\varepsilon^{'}_{n}n^{1/6}]$. Letting $b=x/\sqrt{n}$ and $\tau=(\sqrt{n}/x)^{1/4}$, we have
\begin{eqnarray}\label{eq:j12}
&&J_{1,i} \leq\pr\Big{(}\sum_{k=1}^{n} bY_{k,l}(\theta_{1,i},\theta_{2,i})+\tau b^{2} \widetilde{V}^{2}_{n}\geq x^{2}-\varepsilon_{n}
+\tau(x^{2}+\varepsilon_{n})\Big{)}\cr &&\quad\quad +\pr\Big{(}\sum_{k=1}^{n} -bY_{k,l}(\theta_{1,i},\theta_{2,i})+\tau b^{2}
\widetilde{V}^{2}_{n}\geq x^{2}-\varepsilon_{n} +\tau(x^{2}+\varepsilon_{n})\Big{)}\cr
 &&\quad\leq \pr\Big{(}\sum_{k=1}^{n} b\overline{Y}_{k,l}(\theta_{1,i},\theta_{2,i})+\tau
b^{2} [\widetilde{V}^{2}_{n}-\ep \widetilde{V}^{2}_{n}] \geq x^{2}-\varepsilon_{n} +\tau \widetilde{\varepsilon}_{n}\Big{)}\cr &&\quad\quad+
\pr\Big{(}\sum_{k=1}^{n} -b\overline{Y}_{k,l}(\theta_{1,i},\theta_{2,i})+\tau b^{2} [\widetilde{V}^{2}_{n}-\ep \widetilde{V}^{2}_{n}] \geq
x^{2}-\varepsilon_{n} +\tau \widetilde{\varepsilon}_{n}\Big{)}\cr &&\quad=:J_{3,i},
\end{eqnarray}
where
\begin{eqnarray*}
&&\overline{Y}_{k,l}(\theta_{1,i}=Y_{k,l}(\theta_{1,i},\theta_{2,i})-\ep \overline{Y}_{k,l}(\theta_{1,i},\theta_{2,i}),\cr
&&\widetilde{\varepsilon}_{n}=\varepsilon_{n}+\tau b^{2}(n-\ep \widetilde{V}^{2}_{n})=\varepsilon_{n}+o(1)x^{3}/\sqrt{n}.
\end{eqnarray*}
 Let $\eta_{k}=\overline{Y}_{k,l}(\theta_{1,i}, \theta_{2,i})$ and $\xi_{k}= X^{'2}_{k}-\ep X^{'2}_{k}$. Using $|e^{s}-1-s-s^{2}/2|\leq
|s|^{3}e^{s\vee 0}$, we get
\begin{eqnarray*}
\ep e^{2b\eta_{k}+2\tau b^{2}\xi_{k}}&=&1+2\ep(b\eta_{k}+\tau b^{2}\xi_{k})^{2}+O(1)\ep|b\eta_{k}+\tau b^{2}\xi_{k}|^{3}e^{3}\cr &=&1+2b^{2}\ep
\eta^{2}_{k}+4\tau b^{3}\ep(\eta_{k} \xi_{k})+2\tau^{2}b^{4}\ep\xi^{2}_{k}+O(1)e^{3}(b^{3}\ep|\eta_{k}|^{3}+\tau^{3}b^{6}\ep|\xi_{k}|^{3})\cr &=&
1+2b^{2}(\ep X^{'2}_{k}-(\ep X^{'}_{k})^{2})[\theta^{2}_{1,i}\cos^{2}(kw_{l})+\theta^{2}_{2,i}\sin^{2}(kw_{l})\cr
&&+2\theta_{1,i}\theta_{2,i}\cos(kw_{l})\sin(kw_{l})]+O(1)(1+\tau)b^{3}\cr &=&
1+2b^{2}(\theta^{2}_{1,i}\cos^{2}(kw_{l})+\theta^{2}_{2,i}\sin^{2}(kw_{l}) +2\theta_{1,i}\theta_{2,i}\cos(kw_{l})\sin(kw_{l}))\cr &
&+O(1)(1+\tau)b^{3},
\end{eqnarray*}
for $x\in [4, \varepsilon^{'}_{n}n^{1/6}]$. This, together with  (\ref{eq:fat}), implies that
\begin{eqnarray}\label{eq:j13}
J_{3,i}&\leq& 2\exp\Big{(}-2x^{2}+2\varepsilon_{n}-2\tau \widetilde{\varepsilon}_{n}+nb^{2}+O(1)(1+\tau)x^{3}/\sqrt{n}\Big{)}\cr &\leq&
C\exp\Big{(}-x^{2}-b^{-1/8}\Big{)}.
\end{eqnarray}
Combining (\ref{eq:j10})-(\ref{eq:j13}) gives
\begin{eqnarray}\label{eq:j14}
\pr\Big{(}\widetilde{M}_{n}\geq \sqrt{n}x,\textbf{A}\Big{)}=o(1)e^{-x^{2}}.
\end{eqnarray}
Define
\begin{eqnarray*}
A_{j}=\Big{\{}\Big{|}\sum_{k=1}^{n}X^{'}_{k}e^{ik\omega_{j}}\Big{|}\geq x\sqrt{n(1+\varepsilon_{n}/x^{2})}\Big{\}}, \quad 1\leq j\leq q.
\end{eqnarray*}
We have
\begin{eqnarray*}
\pr\Big{(}\widetilde{M}_{n}\geq x\sqrt{n(1+\varepsilon_{n}/x^{2})}\Big{)}\geq \sum_{j=1}^{q}\pr(A_{j})-\sum_{1\leq i<j\leq q}\pr(A_{i}A_{j}).
\end{eqnarray*}
By Lemma \ref{le4.2} (ii), we have
\begin{eqnarray*}
\pr(A_{i})=(1+o(1))e^{-x^{2}},\quad\pr(A_{i}A_{j})=2^{-1}(1+o(1))e^{-2x^{2}},
\end{eqnarray*}
uniformly in $x\in [4, \varepsilon^{'}_{n}n^{1/6}]$ and $1\leq i,j\leq q$. This shows that
\begin{eqnarray}\label{eq10}
~~\pr\Big{(}\widetilde{M}_{n}\geq x\sqrt{n(1+\varepsilon_{n}/x^{2})}\Big{)}\geq (1+o(1))qe^{-x^{2}}(1-2^{-1}qe^{-x^{2}}).
\end{eqnarray}
It follows from (\ref{eq9}), (\ref{eq:j14}) and (\ref{eq10}) that
\begin{eqnarray}\label{p16}
\pr\Big{(}M_{n}\geq xV_{n}\Big{)}\geq (1+o(1))qe^{-x^{2}}(1-2^{-1}qe^{-x^{2}}).
\end{eqnarray}
uniformly in $x\in [4, \varepsilon^{'}_{n}n^{1/6}]$. Let $x=\sqrt{y+\log q}$. Combining (\ref{p15}) and (\ref{p16}), we have for any fixed $M>0$,
\begin{eqnarray}\label{eq:j15}
\limsup_{n\rightarrow\infty}\sup_{M\leq y\leq \varepsilon^{'}_{n}n^{1/3}}\Big{|}\frac{\pr\Big{(}\frac{M^{2}_{n}}{V^{2}_{n}}-\log q\geq
y\Big{)}}{1-\exp(-\exp(-y))}-1\Big{|}\leq Ce^{-M}.
\end{eqnarray}
For $-\log q\leq y\leq M$, by (\ref{eq:f}), (\ref{f1}) and (\ref{f2}),
\begin{eqnarray}\label{eq:j16}
\limsup_{n\rightarrow\infty}\sup_{-\log q\leq y\leq M}\Big{|}\frac{\pr\Big{(}\frac{M^{2}_{n}}{V^{2}_{n}}-\log q\geq
y\Big{)}}{1-\exp(-\exp(-y))}-1\Big{|}=0.
\end{eqnarray}
This proves Theorem 2.3 by (\ref{eq:j15}) and (\ref{eq:j16}). $\Box$\vspace{0.5cm}

\noindent{\bf Proof of Lemma 3.1.} This lemma follows immediately by Theorem 2.3 and
\begin{eqnarray*}
f_{n}((x+\log q)/q)=\pr\Big{(}\frac{\max_{1\leq j\leq q}I_{n}(\omega_{j})}{q^{-1}\sum_{j=1}^{q}I_{n}(\omega_{j})}-\log q\geq x\Big{)},
\end{eqnarray*}
where $\{X_{k}\}$ are i.i.d. $N(0,1)$ random variables. $\Box$\vspace{0.5cm}

\noindent{\bf Proof of Theorem 3.1.}  Let $\mathcal{C}_{g}=\{P_{g}<\theta/(3G), P^{true}_{g}>\theta/(2G)\}$ and define $F(x)=\exp(-\exp(-x))$. Let
$x_{n}$ satisfy $1-F(x_{n})=\theta/(2.5G)$. So $x_{n}\sim \log G$. Corollary 3.1  yields
\begin{eqnarray}\label{eq7}
\max_{1\leq g\leq G}\Big{|}\frac{1-F(x_{n})}{1-F_{n,g}((x_{n}+\log q)/q)}-1\Big{|}=o(1).
\end{eqnarray}
By (\ref{eq7}) and the definition of $x_{n}$, we can see that on $\mathcal{C}_{g}$, it holds $P^{true}_{g}>\theta/(2G)>1-F_{n,g}((x_{n}+\log
q)/q)$ for $n$ large. By the monotonicity of distribution function we have $qf_{g}-\log q\leq x_{n}$. This together with Corollary 3.1 and Lemma
3.1 yields
\begin{eqnarray*}
\max_{1\leq g\leq G}\Big{|}\frac{P_{g}}{P^{true}_{g}}-1\Big{|}I\{\mathcal{C}_{g}\}=o(1).
\end{eqnarray*}
Note that on $\mathcal{H}_{g}\bigcap \mathcal{C}^{c}_{g}$ we have $P_{g}\geq\theta/(3G)$. We can show  that $qf_{g}-\log q\leq y_{n}$, where
$1-F(y_{n})=\theta/(4G)$, so  $y_{n}\sim \log G$. In fact, by Lemma 3.1,
\begin{eqnarray*}
\frac{f_{n}(q^{-1}(y_{n}+\log q))}{1-F(y_{n})}-1=o(1),
\end{eqnarray*}
and hence, $f_{n}(f_{g})=P_{g}>f_{n}(q^{-1}(y_{n}+\log q))$ for $n$ large, which implies  $qf_{g}-\log q\leq y_{n}$. It follows from Corollary 3.1
and Lemma 3.1 that
\begin{eqnarray*}
\max_{1\leq g\leq G}\Big{|}\frac{P_{g}}{P^{true}_{g}}-1\Big{|}I\{\mathcal{H}_{g}\}I\{\mathcal{C}^{c}_{g}\}=o(1).
\end{eqnarray*}
The theorem is proved. $\Box$


\begin{thebibliography}{99}

\bibitem{ahd} Ahdesm\"{a}ki, M., L\"{a}hdesm\"{a}ki, H., Pearson, R.,
Huttunen, H.  and Yli-Harja, O. (2005). Robust detection of periodic time series measured from biological systems. {\em BMC Bioinformatics,} {\bf
6: 117} 1-18.\vspace{-0.3cm}

\bibitem{amo} Amosova, N.N. (1982). Probabilities of moderate deviations. {\em J. Math. Sci.} {\bf 20:} 2123-2130.\vspace{-0.3cm}

\bibitem{an} An, H. Z., Chen, Z. G. and Hannan, E. J. (1983). The maximum of the periodogram.
{\em J. Multivariate Anal. }{\bf 13: } 383-400.\vspace{-0.3cm}

\bibitem{ben} Benjamini, Y. and Hochberg, Y. (1995). Controlling the
false discovery rate: a practical and powerful approach to multiple testing. {\em J. Roy. Statist. Soc. Ser. B} {\bf 57:} 289-300.\vspace{-0.3cm}







\bibitem{chen} Chen, J. (2005). Identification of significant periodic
genes in microarray gene expression data. {\em BMC Bioinformatics,} {\bf 6: 286} 1-12.\vspace{-0.3cm}

\bibitem{davis} Davis, R. A. and Mikosch, T.
 (1999). The maximum of the periodogram of a non-Gaussian sequence. {\em
 Ann. Probab.} {\bf 27: } 522-536.\vspace{-0.3cm}

\bibitem{dela} de la Pe\"{n}a, V.H., Lai, T.L. and Shao, Q.M. (2009). Self-normalized processes: limit theory and statistical applications.
Springer, New York.\vspace{-0.3cm}


\bibitem{fan} Fan, J.Q., Hall, P and Yao, Q. (2007). To how many
simultaneous hypothesis tests can normal, student's t or bootstrap calibration be applied? {\em J. Am. Stat. Assoc.} {\bf
102:}1282-1288.\vspace{-0.3cm}



\bibitem{fay} Fay, G. and Soulier, P. (2001). The periodogram of an i.i.d.
sequence. {\em Stoch. Proc. Appl.} {\bf 92: } 315-343.\vspace{-0.3cm}

\bibitem{fish} Fisher, R. A. (1929). Tests of significance in harmonic analysis. {\em Proc. Roy. Statist. Soc. Ser. A } {\bf 125:} 54-59.\vspace{-0.3cm}

\bibitem{gly} Glynn, E., Chen, J. and Mushegian, A. (2006). Detecting
periodic patterns in unevenly spaced gene expression time series using Lomb-Scargle periodograms. {\em Bioinformatics,} {\bf 22(3):}
310-316.\vspace{-0.3cm}

\bibitem{hu} Hu, Z., Shao, Q.M. and Wang, Q. (2009).
Cram\'{e}r-type moderate deviations for the maximum of self-normalized sums. {\em Elect. J. Probab.} {\bf 14:} 1181-1197.\vspace{-0.3cm}


\bibitem{jing} Jing, B.Y., Shao, Q.M. and Wang, Q. (2003). Self-normalized
Cram\'{e}r-type large deviations for independent random variables. {\em Ann. Probab.} {\bf 31:} 2167-2215.\vspace{-0.3cm}

\bibitem{korosok} Korosok, M.R. and Ma, S. (2007), Marginal asymptotics for
the large $p$, small $n$ paradigm: with applications to microarray data. {\em Ann. Statist.} {\bf 35:} 1456-1486.\vspace{-0.3cm}


\bibitem{ledoux} Ledoux, M. and Talagrand, M. (1991). Probability
in Banach spaces. Springer, Berlin.\vspace{-0.3cm}



\bibitem{lin1} Lin, Z.Y. and Liu, W.D. (2009a). On maxima of periodograms
of stationary processes. {\em Ann. Statist.} {\bf 37: } 2676-2695.\vspace{-0.3cm}

\bibitem{lin2}  Lin, Z.Y. and Liu, W.D. (2009b). Supplementary material for "On maxima of periodograms of stationary processes."
Available at
 http://xxx. arxiv.org/abs/0801.1357.\vspace{-0.3cm}

\bibitem{mich} Michel, R. (1976). Nonuniform central limit bounds with application to probabilities of deviations.
{\em Ann. Probab.}  {\bf 4:} 102-106.\vspace{-0.3cm}

\bibitem {miko} Mikosch, T.,  Resnick, S., and Samorodnitsky, G. (2000).
The maximum of the periodogram for a heavy-tailed sequence. {\em Ann. Probab.} {\bf 28:} 885-908.\vspace{-0.3cm}


\bibitem {Petrov}  Petrov, V.V. (2002). On probabilities of moderate deviations.
{\em J. Math. Sci.} {\bf 109:} 2189-2191.\vspace{-0.3cm}




\bibitem {shao1} Shao, Q. M. (1997). Self-normalized large deviations. {\em Ann.Probab.} {\bf 25:} 285-328.\vspace{-0.3cm}

\bibitem {shao2} Shao, Q. M. (1999). A Cram\'{e}r type large deviation for Student's $t$ statistic. {\em
J. Theoret. Probab.} {\bf 12: } 385-398.\vspace{-0.3cm}

\bibitem{shaox} Shao, X. and Wu, W. B. (2007). Asymptotic spectral theory for nonlinear time
series. {\em Ann. Statist.} {\bf 35 (4):} 1773-1801.\vspace{-0.3cm}

\bibitem{wang} Wang, Q. (2005). Limit theorems for self-normalized large deviations.
{\em Elect. J. Probab.} {\bf 10:} 1260-1285.\vspace{-0.3cm}

\bibitem{wich}  Wichert, S., Fokianos, K. and Strimmer, K. (2004).
Identifying periodically expressed transcripts in microarray time series data. {\em Bioinformatics,} {\bf 20(1):} 5-20.\vspace{-0.3cm}

\bibitem{wu} Wu, W. B. and Zhao, Z. (2008). Moderate deviations for stationary processes. {\em Statist. Sinica.}
{\bf 18: } 769-782.\vspace{-0.3cm}










\end{thebibliography}
\end{document}